\documentclass{article}

\usepackage[english]{babel}

\usepackage[utf8]{inputenc}
\usepackage[T1]{fontenc}

\usepackage{amssymb}
\usepackage{stmaryrd}
\usepackage{amsmath,amsfonts}
\usepackage[tikz]{bclogo}
\usepackage[top=4.34cm]{geometry}

\usepackage{lmodern}
\usepackage{textcomp}
\usepackage{mathtools, bm}
\usepackage{amssymb, bm}
\usepackage[unq]{unique}
\usepackage{enumerate}
\usepackage{comment}
\usepackage[breaklinks]{hyperref}

\usepackage[numbers]{natbib}  

\usepackage[breaklinks]{hyperref}
\usepackage[numbers]{natbib}

\usepackage{amsthm}

\usepackage{cleveref}

\newtheorem{theorem}{Theorem}[section]
\newtheorem{lemma}[theorem]{Lemma}

\newtheorem{conjecture}[theorem]{Conjecture}

\newtheorem{Question}[theorem]{Question}
\newtheorem{Claim}[theorem]{Claim}

\theoremstyle{definition}

\newcommand{\Addresses}{{
  \bigskip
  \footnotesize

  \textsc{Centre for Mathematical Sciences,
Wilberforce Road,
Cambridge CB3 0WA,
United Kingdom}\par\nopagebreak
  \textit{E-mail address:} \texttt{\{jp895,jp899,lvv23\}@cam.ac.uk}

}}

\DeclareMathOperator{\Dom}{Dom}

\title{On the number of minimum dominating sets and total dominating sets in forests}

\author{Jan Petr \and Julien Portier \and Leo Versteegen}

\date{}

\begin{document}

\maketitle

\begin{abstract}
We show that the maximum number of minimum dominating sets of a forest with domination number $\gamma$ is at most $\sqrt{5}^{\gamma}$ and construct for each $\gamma$ a tree with domination number $\gamma$ that has more than $\frac{2}{5}\sqrt{5}^{\gamma}$ minimum dominating sets. Furthermore, we disprove a conjecture about the number of minimum total dominating sets in forests by Henning, Mohr and Rautenbach.
\end{abstract}

\section{Introduction}

Consider a set of vertices $D \subset V(G)$ in a graph $G$. We call $D$ \emph{dominating} if for every vertex  $v\in V(G)$, there exists $w\in V(G)$ such that either $v=w$ or $w$ is a neighbour of $v$. We say then that $w$ \emph{dominates} $v$, noting that a vertex $v$ may be dominated by several vertices in $G$. More restrictively, we say that $D$ is \emph{total dominating} if every vertex in $G$ has a neighbour in $D$.
The minimum size of a dominating set of $G$ is called the \emph{domination number of $G$} and is denoted $\gamma(G)$. Analogously, the minimum size of a total dominating set of $G$ is called the \emph{total domination number of $G$} and is denoted by $\gamma_t(G)$. A minimum dominating set, i.e., a dominating set with $\gamma(G)$ elements, will also be referred to as \emph{$\gamma$-set}. We write $\Dom(G)$ for the set of all $\gamma$-sets of $G$ and $\Gamma(G)$ for $\vert \Dom(G)\vert$. Similarly, $\Gamma_t(G)$ will stand for the number of minimum total dominating sets of $G$.

Our paper is concerned with minimum dominating and total dominating sets in forests. Fricke, Hedetniemi, Hedetniemi and Hutson \cite{Fricke}, when investigating the structure of minimum dominating sets, asked whether the number of minimum dominating sets of a tree of domination number $\gamma$ was at most $2^\gamma$. Edwards, MacGillivray and Nasserasr \cite{Edwards} and independently Bie\'{n} \cite{Bien} gave a negative answer to this question. Edwards et al. \cite{Edwards} also proved an upper bound of the form $c^\gamma$ on the number of minimum dominating sets in a forest:

\begin{theorem}
A forest of domination number $\gamma$ has at most $(\frac{1}{2}(1+\sqrt{13}))^\gamma$ minimum dominating sets.
\end{theorem}

Independently, Alvarado, Dantas, Mohr and Rautenbach \cite{Alvarado} obtained the same result with $\frac{1}{2}(1+\sqrt{13})\approx 2.3028$ replaced by the largest zero $r\approx 2.4605$ of $x^3-x^2-4x+1$, but a refinement of Claim $3$ from their paper would also result in the bound $(\frac{1}{2}(1+\sqrt{13}))^\gamma$.

In this paper we determine the correct order of the maximum number of minimum dominating sets in a forest with domination number $\gamma$. To this end, we first prove for all forests $F$ the upper bound $\Gamma(F)\leq \sqrt{5}^\gamma$.

\begin{theorem}\label{th:UB}
All forests of domination number $\gamma$ have at most $\sqrt{5}^\gamma$ minimum dominating sets.
\end{theorem}

By means of an explicit construction, we show that this upper bound is best possible up to a multiplicative constant.

\begin{theorem}
\label{LowerBound}
For every positive integer $\gamma$, there exists a tree with domination number $\gamma$ that has more than $\frac{2}{5} \sqrt{5}^{\gamma}$ minimum dominating sets.
\end{theorem}

This result also disproves the following conjecture by Alvarado et al. \cite{Alvarado}:

\begin{conjecture}
A tree with domination number $\gamma$ has $O(\frac{\gamma2^\gamma}{\ln\gamma})$ minimum dominating sets.
\end{conjecture}

Results of similar nature as \Cref{th:UB} and \Cref{LowerBound} were recently proven in the setting of minimal (rather than minimum) dominating sets of a tree by G\"{u}nter Rote \cite{Rote}, who showed that a tree on $n$ vertices has at most $\sqrt[13]{95}^n$ minimal dominating sets and that $\sqrt[13]{95}$ is the best possible growth constant. For the case of a general graph on $n$ vertices, Fomin, Grandoni, Pyatkin and Stepanov \cite{Fomin} showed that the number of minimal dominating sets is at most $1.7159^n$ and can be as high as ${15}^{n/6}\approx 1.5704^n$, as given by Kratsch's example of disjoint copies of octahedrons.

In the second part of our paper, we turn our attention to the maximum number of minimum total dominating sets of forests (with no isolated vertices) on $n$ vertices with total domination number $\gamma_t$. This topic was previously analyzed by Henning, Mohr and Rautenbach \cite{Henning}. Here, e.g. in view of stars, an upper bound only in terms on $\gamma_t$ cannot exist. Henning at al. proved the following upper bounds:

\begin{theorem}\label{th:total_ub}
Every forest $F$ with
order $n$, no isolated vertex, and total domination number $\gamma_t$ has at most
$$\min\{(8 \sqrt{e})^{\gamma_t} \left(\frac{n-\frac{\gamma_t}{2}}{\frac{\gamma_t}{2}} \right)^\frac{\gamma_t}{2}, (1+\sqrt{2})^{n-\gamma_t}, 1.4865^n\}$$
minimum total dominating sets.
\end{theorem}

Conjecturing that they found extremal examples for given $n$ and $\gamma_t$, Henning et al.\ also proposed the following upper bound:

\begin{conjecture}\label{Conj_Total_UB}
If a tree $T$ has order $n \geq 2$ and total domination number $\gamma_t$, then $$\Gamma_t(T) \leq \left(\frac{n-\frac{\gamma_t}{2}}{\frac{\gamma_t}{2}} \right)^\frac{\gamma_t}{2}.$$
\end{conjecture}

Our last result disproves this conjecture:

\begin{theorem}\label{th:total_LB}
There exists $c>1$ and trees $T$ with arbitrarily large total domination number $\gamma_t$ and order $n$ with
$$\Gamma_t(T) \geq c^{\gamma_t} \left(\frac{n-\frac{\gamma_t}{2}}{\frac{\gamma_t}{2}} \right)^\frac{\gamma_t}{2}.$$
\end{theorem}

\section{Minimum dominating sets in forests}

Our proof of \Cref{th:UB} uses a general approach similar to the one used by Alvarado et al.\ \cite{Alvarado} to prove their Theorem $2$. Claims \ref{claim:strong-support}, \ref{x_not_support} and \ref{y_not_support_vertex} and their proofs are very close to Claims $1$, $4$ and $6$ and their proofs in \cite{Alvarado}. As there are differences in our notations, we prove all claims in full.

Following Alvarado et al., we will call a vertex $v$ of a forest a \emph{support vertex} if it is a neighbour of a leaf. Furthermore, we call $v$ a \emph{strong support vertex} if it is adjacent to at least two leaves. Note that every strong support vertex is contained in any minimum dominating set.

Moreover, we define inductively a $k$-terminal vertex of a tree. We say that a leaf is \emph{$0$-terminal} and for all positive integers $k$ that a vertex $v$ is \emph{$k$-terminal} if at most one of its neighbours is not $l$-terminal for some $l < k$ and $v$ is not $l$-terminal for some $l<k$. Let $L(F)$ stand for the set of leaves of a forest $F$. One can equivalently define $k$-terminal vertices as follows. Starting with $F_0=F$, for all non-negative integers $k$ define $F_{k+1}=F-L(F_{k})$. The set of all $k$-terminal vertices is then precisely $L(F_k)$. Note that for any $k \geq 1$ a $k$-terminal vertex has at least two neighbours, out of which at least one is $(k-1)$-terminal.

By a \emph{branch} of a forest $F$ from a vertex $v$ we mean a component $C$ of $F-v$ which contains a neighbour of $v$ in $F$. For a branch $B$ from a vertex $v$ we define its \emph{depth} as the maximum distance between $v$ and a vertex from $B$. A $k$-terminal vertex $v$ thus always has at most one branch from $v$ of depth greater than $k$.

\subsection{Proof of the upper bound}

For a positive integer $\gamma$ and a non-negative integer $s$ satisfying $s \leq \gamma$, let $f(\gamma, s)$ be the maximum number of minimum dominating sets among all forests with domination number $\gamma$ and at least $s$ strong support vertices.\footnote{We remark that the original definition in \cite{Alvarado} did not include the term `at least'. We have made this change because it is convenient to have $f$ be non-increasing in $s$.} Note that $f(\gamma, s)$ is non-increasing in $s$.

Let $\beta = \sqrt{5}$ and $\alpha=\frac{\beta}{\beta -1}$. We will prove that $f(\gamma,s) \leq \alpha^{s} \beta^{\gamma -s}$. For a contradiction, let $\gamma$ and $s$ be such that $f (\gamma, s) > \alpha^{s} \beta^{\gamma -s}$ and $\gamma$ is as small as possible. Let $T$ be a forest with domination number $\gamma$ and at least $s$ strong support vertices that has $f (\gamma, s)$ minimum dominating sets. Note that $T$ is actually a tree, as otherwise there would exist a component of $T$ with domination number $\gamma'<\gamma$ and $s'\leq s$ support vertices such that $f (\gamma', s') > \alpha^{s'} \beta^{\gamma' -s'}$, contradicting the minimality of $\gamma$. Furthermore, $\gamma>1$, as the only tree with domination number $1$ that has more than $1$ minimum dominating set is a path on two vertices.

From here on out we will lead the proof as follows. We first prove a series of claims showing that $T$ has to have a simple, linear structure near some of its leaves. Each of the claims will assume the existence of a branching at one or more vertices at distance up to $5$ from a leaf. This assumption on the structure of $T$, together with the fact that $f(\gamma',s')\leq \alpha^{s'}\beta^{\gamma'-s'}$ for all $\gamma'<\gamma$, will allow us to derive an upper bound on $\Gamma(T)$, which will be in contradiction to the lower bound $\Gamma(T)>\alpha^s \beta^{\gamma-s}$. Once we are equipped with enough structural information on $T$, we can give an assumption-free upper bound on $\Gamma(T)$, which will nonetheless contradict the lower bound $\Gamma(T)>\alpha^s \beta^{\gamma-s}$, thus completing the proof.

\begin{Claim}\label{claim:strong-support}
A $1$-terminal vertex has degree $2$, i.e., it is not a strong support vertex.
\end{Claim}

\begin{proof}
We assume for contradiction that there exists a $1$-terminal strong support vertex $v$. As $\gamma$ is at least $2$, $v$ has a non-leaf neighbour $w$ (and not more than one as $v$ is $1$-terminal), we distinguish whether $w$ is a strong support vertex or not. If it is, $w$ and $v$ must be in $D$, and therefore any element $D$ of $\Dom(T)$ has the form $D'\cup \{v\}$ for a set $D'\in \Dom( T-(\{v\}\cup N(v)\setminus \{w\}))$ of size $\gamma-1$. There are at most $f(\gamma-1,s-1)$ such sets $D'$ which is at most $\alpha^{s-1}\beta^{\gamma-s}<\alpha^{s}\beta^{\gamma-s}$ by induction. This contradicts our assumed lower bound on $\Gamma(T)$.

If $w$ is not a strong support vertex, then we split the sets in $\Dom(T)$ into two types, those that contain $w$ and those that do not. The latter are always of the form $D_1\cup \{v\}$ for some $D_1\in \Dom(T-(\{v\} \cup N(v)))$ of size $\gamma-1$ of which there are at most $f(\gamma-1,s-1)$ many. The former are of the form $D_2\cup \{v\}$ for some $D_2\in \Dom(T-(\{v\} \cup N(v)\setminus\{w\}))$ of size $\gamma-1$ that contains $w$. Since $w$ is not a strong support vertex of $T$, the number of such sets $D_2$ is given by the number of (minimum) dominating sets of size $\gamma-1$ in the tree $T'$ that is obtained from $T$ by removing $v$ and its neighbouring leaves and appending two leaves to $w$, so that $w$ becomes a strong support vertex. In total, we have at most $f(\gamma-1,s-1)+f(\gamma-1,s)\leq \alpha^s\beta^{\gamma-s}(\alpha^{-1}+\beta^{-1})=\alpha^s\beta^{\gamma-s}$ minimum dominating sets for $T$, again contradicting our assumption on $\Gamma(T)$.
\end{proof}

\begin{Claim}\label{Only_forks_or_3paths}
Let $w$ be a $2$-terminal vertex. Then all but one branches from $w$ must be paths of length $2$ and $w$ can have at most $3$ neighbours. Moreover, $w$ has a $l$-terminal neighbour for $l \geq 2$.
\end{Claim}

\begin{figure}[htbp]\centering
    			\includegraphics[height=3cm]{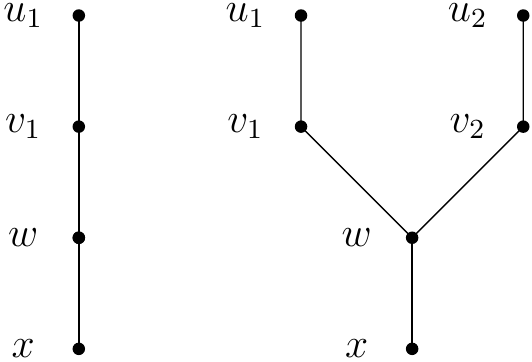}
    			\caption{Illustration of \Cref{Only_forks_or_3paths}: each $2$-terminal vertex $w$ has one or two $1$-terminal neighbours ($v_1,v_2$), each with one leaf neighbour ($u_1,u_2$), and another neighbour $x$ which is $l$-terminal for $l \geq 2$.}
\label{fig:claim2.2}
	\end{figure}

\begin{proof}
As $w$ is $2$-terminal, it has a $1$-terminal neighbour $v$, which by \Cref{claim:strong-support}, has exactly one leaf $u$. 
Suppose now that $w$ has a leaf neighbour $v'$ as well. Any minimum dominating set $D$ of $T$ must contain exactly one vertex in $\{u,v\}$, call it $r$. It also must contain exactly one vertex in $\{w,v'\}$, call it $q$. We claim that $D \setminus \{r\}\in \Dom(T-\{u,v\})$. Indeed, the set $D \setminus \{r\}$ is dominating for $T-\{u,v\}$, as $w$, the only vertex in $T-\{u,v\}$ that could be a neighbour of $r$, is dominated by $q$. Furthermore, the set $D \setminus \{r\}$ is a minimum dominating set of $T-\{u,v\}$, as every smaller dominating set of $T-\{u,v\}$ would give rise to a smaller dominating set of $T$ by adding $r$ again. Therefore, we have $\Gamma(T)\leq 2\Gamma(T-\{u,v\})\leq 2 f(\gamma-1,s)$. By induction, the latter is at most $2\alpha^s\beta^{\gamma-1-s} \leq \alpha^s\beta^{\gamma-s}$ as $\beta \geq 2$, which gives a contradiction.

Thus $w$ does not have any neighbouring leaf and we can apply \Cref{claim:strong-support} to every $1$-terminal neighbour of $w$ to see that all branches of $w$ of depth at most $2$ are in fact paths of length $2$. Assume $w$ has $m$ such branches with vertices $u_1,v_1,\ldots,u_m,v_m$ such that $u_iv_iw$ are the respective paths. It remains to show that $m\leq 2$ and $w$ has an $l$-terminal neighbour for $l \geq 2$.

To see this, assume first that $w$ has no other neighbour than the $v_i$'s. Then $\gamma=m$, $s=0$ and all minimum sets are given by including from each branch either $u_i$ or $v_i$ with the requirement that $v_i$ must be chosen at least once, giving $\Gamma(T)=2^{m}-1<\beta^{m}$.

We are left with the case that $w$ has an additional neighbour $x$ which must be $l$-terminal neighbour for $l \geq 2$. Still, each minimum dominating set must contain one element $r_i\in \{u_i,v_i\}$ for every $i\in [m]$. Furthermore, every minimum dominating set $D$ of $T$ satisfies at least one of the following conditions:
\begin{itemize}
    \item $r_i=u_i$ for every $i\in [m]$. In this case, $D\setminus \{u_1,\ldots,u_m\}$ is a minimum dominating set of $T-\{u_1,v_1,\ldots,u_m,v_m\}$.
    \item $r_i=v_i$ for at least one $i\in [m]$ and $D$ does not contain $w$. In this case $D\setminus \{r_1,\ldots,r_m\}$ is a minimum dominating set of $T-\{w,u_1,v_1,\ldots,u_m,v_m\}$.
    \item $r_i=v_i$ for at least one $i\in [m]$ and $D$ does contain $w$. In this case $x$ is not an element of $D$ as otherwise $w$ would be redundant and so $D\setminus \{w,r_1,\ldots,r_m\}$ must be a minimum dominating set of $T-\{x,w,u_1,v_1,\ldots,u_m,v_m\}$. 
\end{itemize}

We see that there are at most $f(\gamma-m,s)$ sets satisfying the first condition, $(2^m-1)f(\gamma-m,s)$ satisfying the second condition and $(2^m-1)f(\gamma-m-1,s)$ satisfying the third condition. By induction, this means that

\begin{align*}
    \Gamma(T)\leq \alpha^s\beta^{\gamma-s} (2^m\beta^{-m}+(2^m-1)\beta^{-m-1}),
\end{align*}
which is less than $\alpha^s\beta^{\gamma-s}$ if $m\geq 3$.
\end{proof}

\begin{Claim}\label{fork_then_x_degree_2}
Let $x$ be a $3$-terminal vertex. If one of its $2$-terminal neighbours has degree more than $2$, then $x$ has degree $2$.
\end{Claim}

\begin{proof}
Suppose a $2$-terminal neighbour $w$ of $x$ has degree more than $2$. By \Cref{Only_forks_or_3paths}, $w$ has exactly two neighbours $v,v'$ other than $x$, both of which are $1$-terminal by \Cref{Only_forks_or_3paths} and have exactly one more neighbour, the leaves $u,u'$ respectively.

Assume for contradiction that $x$ has at least two neighbours other than $w$. If $x$ has an $l$-terminal neighbour for $l\geq 3$, call it $y$. All the neighbours of $x$ except for $w$ and $y$, call them $c_1\ldots,c_m$ $(m \geq 1)$, are $l$-terminal for some $l<3$. First, suppose any of $c_1,\ldots,c_m$ is $0$-terminal or $1$-terminal. Then it is straightforward to see that $w$ cannot be in any minimum dominating set of $T$. Any minimum dominating set of $T$ is thus of the form $D_1 \cup \{u\}$ where $D_1 \in \Dom(T - \{u,v\} )$ or $D_2 \cup \{v\}$ where $D_2 \in \Dom(T-\{u,v,w\})$. We then have $\Gamma(T) \leq 2f(\gamma-1, s) \leq \alpha^{s}\beta^{\gamma -s}$ which yields a contradiction.

Hence, we may suppose that $c_1,\ldots,c_m$ are all $2$-terminal. By \Cref{Only_forks_or_3paths}, each of them has one or two other neighbours, all of which are $1$-terminal vertices of degree $2$. We now claim that if $x$ is in a minimum dominating set $D$ of $T$, then none of the neighbours of $x$ is in $D$. To see this note first that for each $1$-terminal vertex $t$ each dominating set contains either $t$ or its leaf-neighbour. Therefore, if $x$ is in $D$, each of $w, c_1, \ldots, c_m$ and all of their neighbours are dominated by vertices in $D \setminus \{w, c_1, \ldots, c_m\}$, and therefore none of $w, c_1, \ldots, c_m$ is in $D$. Assume now that $y$ is in $D$. Denote by $B$ the union of $\{x,y\}$ with all the vertices of the branches from $x$ that do not contain $y$. Then for each $1$-terminal vertex $t$ among those in $B$, $D \cap B$ contains either $t$, or its leaf neighbour. As it is possible to dominate all of $B$ by all such $1$-terminal vertices and $y$, $x$ could not be in $D$ by the minimality of $|D|$.
We have thus proved that if $x$ is an element of $D$ then none of its neighbours is.

Observe that if $y$ exists, it cannot be a strong support vertex. Indeed, if $y$ were a strong support vertex, $w$ again could not be in any minimum dominating set and we would reach contradiction as before. This information about $y$ will later enable us to remove it without decreasing the number of strong support vertices. The same holds for the $c_i$'s.

Let the neighbours of $c_1$ other than $x$ be $b$ and $b'$ (if it exists) and their respective leaf neighbours be $a$ and $a'$ (if $b'$ exists). We distinguish two cases based on whether $b'$ exists or not (see \Cref{fig:fork_x_one_child}).

\begin{figure}[htbp]\centering
    			\includegraphics[height=3cm]{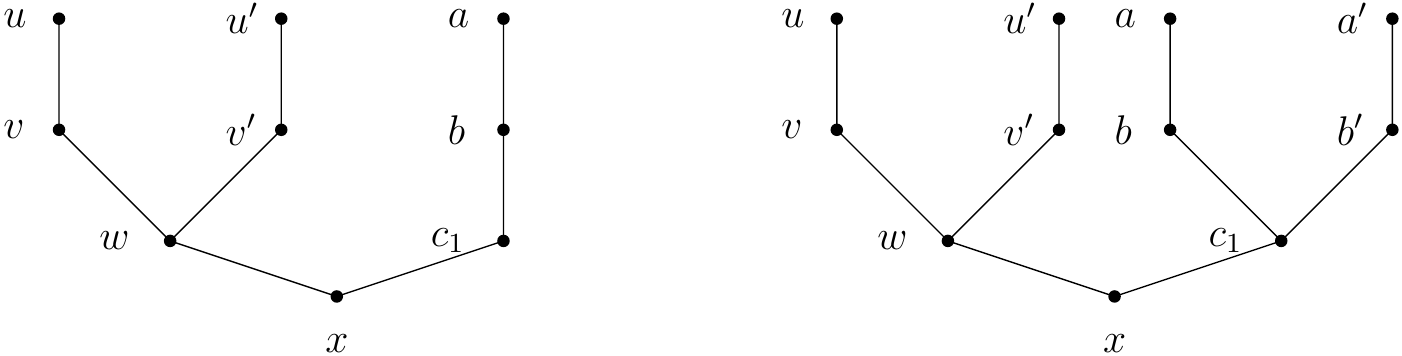}
    			\caption{Proof of \Cref{fork_then_x_degree_2}. In the case that $c_1$ is a $2$-terminal vertex, it has either degree $2$ (on the left), or $3$ (on the right).}
\label{fig:fork_x_one_child}
	\end{figure}

First, suppose $b'$ does not exist. We have already established that no neighbour of $x$ can be contained in $D$ if $D$ contains $x$ itself. Note furthermore that $w$ and $c_1$ cannot both belong to the same minimum dominating set. We may therefore distinguish between the following four possibilities: 

\begin{itemize}
    \item If $c_1\in D$, then $D$ is of the form $D_1 \cup E_1 \cup E_2$ where $D_1 \in \Dom(T - \{ u,u',v,v',w,x,a,b,c_1\})$ of size $\gamma - 4$ (if $\gamma(T - \{ u,u',v,v',w,x,a,b,c_1\})$ is larger than $\gamma-4$, no minimum dominating set of $T$ has this form), $E_1$ is an element of $ \{ \{a,c_1\} ,  \{b,c_1\} \}$ and $E_2$ is an element of $ \{ \{u,v'\} ,  \{u',v\}, \{v,v'\} \}$.
    
    \item If $x\in D$, $D$ is of the form $D_2 \cup F_1 \cup F_2$ where $D_2 \in \Dom(T - (\{ u,u',v,v',x,a,b\} \cup N(x)))$ of size $\gamma - 4$, $F_1$ is an element of $ \{ \{a,x\} ,  \{b,x\} \}$ and $F_2$ is an element of $ \{ \{u,v'\} ,  \{u',v\}, \{v,v'\}, \{u,u'\} \}$.
    
    \item If $w\in D$, $D$ is of the form $D_1 \cup \{ b,w \} \cup F_2$ where $D_1 \in \Dom(T - \{ u,u',v,v',w,x,a,b,c_1\})$ of size $\gamma - 4$ and $F_2$ is as above.
    
    \item If none of $c_1,x,w$ is contained in $D$, $D$ is of the form $D_3 \cup \{ b \} \cup E_2$ where $D_3 \in \Dom(T - \{ u,u',v,v',w,a,b,c_1\})$ of size $\gamma - 3$ and $E_2$ is as above.
    
\end{itemize}

This gives
\begin{align*}
    \Gamma(T)  &\leq 6 f(\gamma-4,s) + 8f(\gamma-4,s)+ 4f( \gamma-4,s) + 3 f(\gamma-3,s)   \\
    &\leq \alpha^{s}\beta^{\gamma -s}( 18\beta^{-4} + 3\beta^{-3}) \\
    &\leq \alpha^{s}\beta^{\gamma -s},
\end{align*}
which gives a contradiction.\\

Finally, suppose $b'$ exists. First, note that $w$ and $c_1$ cannot both belong to a minimum dominating set. A set $D$ is a minimum dominating set of $T$ only if one of the following holds:

\begin{itemize}
    \item either $D$ is of the form $D_1 \cup E_1 \cup E_2$ where $D_1 \in \Dom(T - \{ u,u',v,v',a,a',b,b',w,c_1,x\})$ of size $\gamma - 5$, $E_1$ is an element of $\{c_1\} \times \{ a,b \} \times \{ a',b'\} $ and $E_2$ is an element of $\{ \{u,v' \}, \{v,u' \}, \{v,v'\} \}$,
    
    \item or $D$ is of the form $D_1 \cup F_1 \cup F_2$ where $D_1 \in \Dom (T - \{ u,u',v,v',a,a',b,b',w,c_1,x\})$ of size $\gamma - 5$, $F_1$ is an element of $\{w\} \times \{ u,v \} \times \{ u',v'\} $ and $F_2$ is an element of $\{ \{a,b' \}, \{b,a' \}, \{b,b'\} \}$,
    
    \item or $D$ is of the form $D_2 \cup E_2 \cup F_2$ where $D_2 \in \Dom(T - \{ u,u',v,v',a,a',b,b',w,c_1\})$ of size $\gamma - 4$ and $E_2$ and $F_2$ are as above,
    
    \item or $D$ is of the form $D_3 \cup \{x, u, u', a ,a'\}$ or $D_3 \cup \{x, u, u'\} \cup F_2$ or $D_3 \cup \{x, a, a'\} \cup E_2$ where $D_3 \in \Dom(T - (\{ x, u,u',v,v',a,a',b,b' \} \cup N(x)))$ of size $\gamma - 5$ and $E_2$ and $F_2$ are as above.
\end{itemize}

This gives
\begin{align*}
    \Gamma(T)  &\leq 12 f(\gamma-5,s) + 12 f(\gamma-5,s)+ 9 f(\gamma-4,s)+ 7 f(\gamma-5,s)  \\ 
    &\leq \alpha^{s}\beta^{\gamma -s}(31 \beta^{-5} + 9 \beta^{-4}) \\
    &\leq \alpha^{s}\beta^{\gamma -s},
\end{align*}
which yields a contradiction, and that finishes the proof of this claim.
\end{proof}

\begin{Claim}\label{x_not_support}
Let $x$ be a $3$-terminal vertex. If one of its $2$-terminal neighbours has degree $2$, then $x$ is not a support vertex.
\end{Claim}

\begin{proof}
Let $w$ be a $2$-terminal neighbour of $x$ of degree $2$ and suppose for contradiction that $x$ is a support vertex, i.e., has a leaf $w'$ adjacent to it. 

Note first that the vertex $x$ must have at least one further neighbour $y$ as otherwise $x$ would be $1$-terminal. By \Cref{fork_then_x_degree_2}, $w$ has therefore only one further neighbour $v$, which is $1$-terminal, and by \Cref{claim:strong-support}, $v$ has only one neighbouring leaf $u$.

First observe that $w$ can never be in such a set. Indeed, any $D\in \Dom(T)$ must contain both an element $r_1\in \{a,x\}$ and $r_2\in \{u,v\}$ and the set $\{v,x\}$ dominates everything that $\{r_1,r_2,w\}$ dominates while using one less vertex. Therefore, any minimum dominating set of $T$ is of the form $D_1 \cup \{u\}$ where $D_1 \in \Dom(T - \{u,v\} )$ of size $\gamma - 1$ or $D_2 \cup \{v\}$ where $D_2 \in \Dom(T-\{u,v,w\})$ of size $\gamma - 1$. Overall, $\Dom(T)$ can therefore have at most $2f(\gamma-1,s) \leq \alpha^s\beta^{\gamma-s}$ elements, yielding a contradiction.
\end{proof}

\begin{Claim}\label{path_x_deg_2}
Every $3$-terminal vertex has degree $2$.
\end{Claim}

\begin{proof}
Let $x$ be a $3$-terminal vertex. By \Cref{fork_then_x_degree_2}, it is enough to show that if one of the $2$-terminal neighbours of $x$ has degree $2$, then $x$ has degree $2$.
Let $w$ be a $2$-terminal neighbour of $x$ that has degree $2$ and suppose, for a contradiction, that $x$ has degree at least $3$. Then $x$ has a neighbour $w'$ different from $w$ that is $k$-terminal with $k \leq 2$. By \Cref{x_not_support}, we know that $k \geq 1$. First, suppose $k=2$, and let $v$ and $v'$ be $1$-terminal neighbours of $w$ and $w'$ respectively. By \Cref{fork_then_x_degree_2}, neither $w$ nor $w'$ can have further neighbours, and by \Cref{claim:strong-support}, $v$ and $v'$ have one leaf each which we call $u$ and $u'$ (see \Cref{fig:path_x_deg_2}).

\begin{figure}[htbp]\centering
    			\includegraphics[height=3cm]{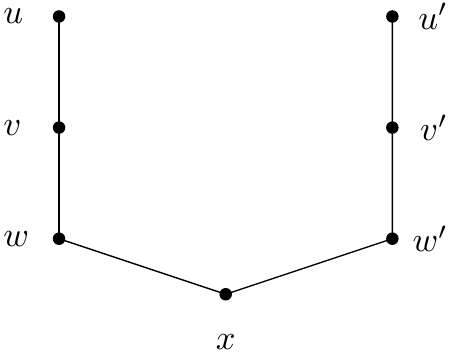}
    			\caption{Proof of \Cref{path_x_deg_2}. Case $k=2$.}
\label{fig:path_x_deg_2}
	\end{figure}

Note that each $D \in \Dom(T)$ contains at most one element from $x$, $w$ and $w'$, as otherwise $D'=(D \setminus \{x,w,w'\}) \cup \{x\})$ would be a smaller minimum dominating set of $T$. This shows that a set $D$ is a minimum dominating set of $T$ if and only if one of the following holds:

\begin{itemize}
\item either $D$ is of the form $D_1 \cup \{ u,u' \}$ where $D_1$ is a minimum dominating set of $T - \{ u,v,u',v' \}$ of size $\gamma - 2$,

\item or $D$ is of the form $D_2 \cup \{ u,v' \}$ where $D_2$ is a minimum dominating set of $T - \{ u,v,u',v',w' \}$ of size $\gamma - 2$,

\item or $D$ is of the form $D_3 \cup \{ v,u' \}$ where $D_3$ is a minimum dominating set of $T - \{ u,v,u',v',w \}$ of size $\gamma - 2$,

\item or $D$ is of the form $D_4 \cup \{ v,v' \}$ where $D_4$ is a minimum dominating set of $T - \{ u,v,u',v',w,w' \}$ of size $\gamma - 2$,

\item or $D$ is of the form $D_5 \cup \{ v,v',w \}$ or $D_5 \cup \{ v,v',w' \}$ where $D_5$ is a minimum dominating set of $T - \{ u,v,u',v',w,w',x \}$ of size $\gamma - 3$.

\end{itemize}

Therefore we get

\begin{align*}
    \Gamma(T)  &\leq f(\gamma-2,s+1) + f(\gamma-2,s) + f(\gamma-2,s) + f(\gamma-2,s) +2f(\gamma-3,s)  \\
    &\leq \alpha^{s}\beta^{\gamma -s}( \beta^{-2}(\beta-1)^{-1}+3\beta^{-2} + 2\beta^{-3}) \\
    &\leq \alpha^{s}\beta^{\gamma -s},
\end{align*}
which gives a contradiction. \\

Now, if $k=1$, by \Cref{claim:strong-support} we know that $w'$ has one leaf $v'$. Note that $w$ cannot belong to a minimum dominating set of $T$, as all of $\{x,w,v,u,w',v'\}$ can be dominated by only $v$ and $w'$ whereas any dominating set containing $w$ would have to intersect $\{w,v,u,w',v'\}$ in at least three vertices. Distinguishing whether $u$ or $v$ belong to the minimum dominating set of $T$, we get $\Gamma(T) \leq 2f( \gamma-1, s) \leq \alpha^{s}\beta^{\gamma -s}$, which gives a contradiction.
\end{proof}

\begin{Claim}\label{fork_then_y_degree_2}
Let $y$ be a $4$-terminal vertex, $x$ be a $3$-terminal neighbour of $y$ and $w$ be a $2$-terminal neighbour of $x$ having degree more than $2$. Then $y$ has degree $2$.
\end{Claim}

\begin{proof}
Suppose for a contradiction that $y$ has a $k$-terminal neighbour $x'$ with $k \leq 3$ different from $x$. If $k=0$, then note that $w$ cannot be in a minimum dominating set of $T$, hence by distinguishing whether $u$ or $v$ belongs to the minimum dominating set, we get that $\Gamma(T)  \leq 2f(\gamma-1, s) \leq \alpha^{s}\beta^{\gamma -s}$, which gives a contradiction. 

\begin{figure}[htbp]\centering
    			\includegraphics[height=4cm]{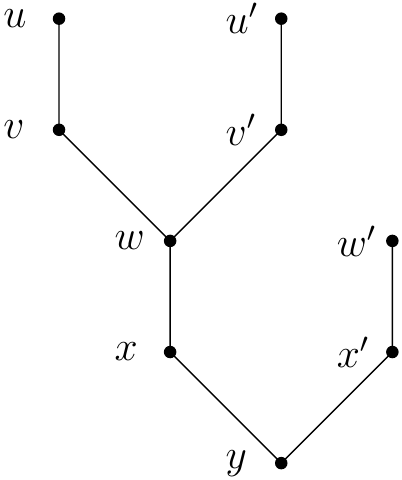}
    			\caption{Proof of \Cref{fork_then_y_degree_2}. Case $k=1$.}
\label{fig:fork_y_one_child}
	\end{figure}

If $k=1$, then by \Cref{claim:strong-support} we know that $x'$ has exactly one $0$-terminal neighbour $w'$ (see \Cref{fig:fork_y_one_child}). We may assume that $y$ is not a strong support vertex as we ruled out the case $k=0$ above. Hence a set $D$ is a minimum dominating set of $T$ if and only if one of the following holds:

\begin{itemize}
\item either $D$ is of the form $D_1 \cup \{ w' \}$ where $D_1$ is a minimum dominating set of $T - \{ w',x' \}$ of size $\gamma - 1$,

\item or $D$ is of the form $D_2 \cup \{ x',u,v' \}$, $D_2 \cup \{ x',v,u' \}$ or $D_2 \cup \{ x',v,v' \}$ where $D_2$ is a minimum dominating set of $T'$ of size $\gamma - 3$, where $T'$ is the tree obtained from $T$ by deleting the set of vertices $\{ w',x',u,v,u',v',w,x \}$ and adding two leaves to $y$,

\item or $D$ is of the form $D_3 \cup \{ x' \} \cup E_1 \cup E_2 \cup E_3$ where $D_3$ is a minimum dominating set of $T - \{ w',x',u,v,u',v',w,x,y \}$ of size $\gamma - 4$, $E_1$ is an element of $\{ \{ u\},\{ v\} \}$,  $E_2$ is an element of $\{ \{ u'\},\{ v'\} \}$, $E_3$ is an element of $\{ \{ w\},\{ x\} \}$.

\end{itemize}

This gives:

\begin{align*}
    \Gamma(T)  &\leq f(\gamma-1,s) + 3f(\gamma-3,s+1) + 8 f(\gamma-4,s)  \\
    &\leq \alpha^{s}\beta^{\gamma -s}( \beta^{-1} +3\beta^{-3}(\beta-1)^{-1} + 8\beta^{-4}) \\
    &\leq \alpha^{s}\beta^{\gamma -s},
\end{align*}
which gives a contradiction. \\

If $k=2$, then let $w'$ be a $1$-terminal neighbour of $x'$, $v''$ be a neighbouring leaf of $w'$. Then note that $x'$ cannot be in a minimum dominating set of $T$. Hence, once again distinguishing whether $v''$ or $w'$ is in a minimum dominating set of $T$, we get a contradiction. \\

Finally, if $k=3$, then by Claims \ref{claim:strong-support}, \ref{Only_forks_or_3paths} and \ref{path_x_deg_2} we know that
$x'$, the $3$-terminal neighbour of $y$ different than $x$, has exactly one $2$-terminal neighbour $w'$, whose one or two other neighbours are $1$-terminal vertices of degree $2$. That forces $y$ to be an element of any minimum dominating set of $T$. Now a minimum dominating set $D$ of $T$ must be of the form $D_1 \cup \{ u,v' \}$, $D_1 \cup \{ v,v' \}$ or $D_1 \cup \{ v,u' \}$ where $D_1$ is a minimum dominating set of $T - \{ u,v,u',v',w \}$ of size $\gamma - 2$. This gives $\Gamma(T) \leq 3f( \gamma-2, s) \leq \alpha^{s}\beta^{\gamma -s}$, which gives a contradiction.
\end{proof}

\begin{Claim}\label{2_term_has_deg_2}
Let $y$ be a 4-terminal vertex, $x$ a 3-terminal neighbour of $y$ and $w$ a $2$-terminal neighbour of $x$. Then $w$ has degree 2.
\end{Claim}

\begin{proof}
Suppose for contradiction that $w$ has degree at least $3$. By Claim \ref{fork_then_y_degree_2}, $y$ has exactly one neighbour $z$ other than $x$ and by Claim \ref{path_x_deg_2}, $x$ has no other neighbours than $y$ and $w$.

\begin{figure}[htbp]\centering
    			\includegraphics[height=4cm]{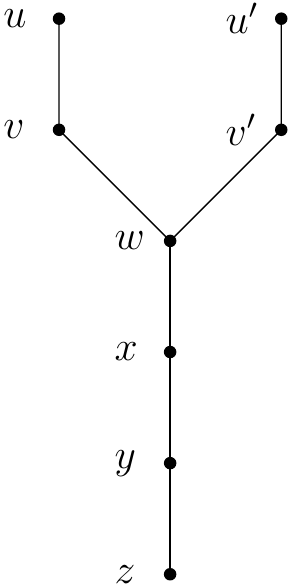}
    			\caption{Proof of \Cref{2_term_has_deg_2}}
\label{fig:fork_z}
	\end{figure}

Suppose that $z$ is not a strong support vertex. Then a set $D$ is a minimum dominating set of $T$ if and only if one of the following holds:

\begin{itemize}
    \item either $D$ is of the form $D_1 \cup \{ u,u',w\}$, or $D_1 \cup \{ u,v',w\}$, or $D_1 \cup \{ v,u',w\}$, or $D_1 \cup \{ v,v',w\}$ where $D_1$ is a minimum dominating set of $T'$ of size $\gamma - 3$ where $T'$ is the tree obtained from $T$ by deleting the set of vertices $\{ u,u',v,v',w,x \}$ and adding two leaves to $z$,
    
    \item or $D$ is of the form $D_2 \cup \{ u,u',x\}$ where $D_2$ is a minimum dominating set of \linebreak $T-\{ u,u',v,v',w,x,y \}$ of size $\gamma - 3$,
    
    \item or $D$ is of the form $D_3 \cup \{ u,v' \}$, or $D_3 \cup \{ u',v\}$, or $D_3 \cup \{ v,v' \}$ where $D_3$ is a minimum dominating set of $T - \{ u,u',v,v',w \}$ of size $\gamma - 2$.
    
\end{itemize}
    
This gives:

\begin{align*}
    \Gamma(T)  &\leq 4 f(\gamma-3,s+1) + f(\gamma-3,s) + 3 f(\gamma-2,s)  \\
    &\leq \alpha^{s}\beta^{\gamma -s}(4 \beta^{-3}(\beta-1)^{-1} + \beta^{-3} + 3 \beta^{-2}) \\
    &\leq \alpha^{s}\beta^{\gamma -s},
\end{align*}
which gives a contradiction. \\

Now, if $z$ is a strong support vertex, then a set $D$ is a minimum dominating set of $T$ if and only if it is of the form $D_4 \cup \{y,v,u' \}$, $D_4 \cup \{y,u,v' \}$, $D_4 \cup \{y,v,v' \}$, $D_4 \cup \{x,u,u' \}$, $D_4 \cup \{x,u,v' \}$, $D_4 \cup \{x,v,u' \}$, $D_4 \cup \{x,v,v' \}$, $D_4 \cup \{w,u,u' \}$, $D_4 \cup \{w,u,v' \}$, $D_4 \cup \{w,v,u' \}$ or $D_4 \cup \{w,v,v' \}$ where $D_4$ is a minimum dominating set of $T - \{ u,u',v,v',w,x,y \}$ of size $\gamma - 3$. This gives $f( \gamma, s) \leq 11f(\gamma-3,s) \leq \alpha^{s}\beta^{\gamma -s}$ which gives a contradiction in this case too.
\end{proof}

\begin{Claim}\label{y_not_support_vertex}
No $4$-terminal vertex $y$ is a support vertex.
\end{Claim}

\begin{proof}
Let $x$ be a $3$-terminal neighbour of $y$, $w$ a $2$-terminal neighbour of $x$ and $v$ a $1$-terminal neighbour of $w$. By \Cref{2_term_has_deg_2}, $w$ has degree $2$ and by \Cref{path_x_deg_2}, so does $x$. By \Cref{claim:strong-support}, $v$ also has only one neighbour other than $w$, a leaf $u$.

Suppose for a contradiction that $y$ is a support vertex and let $x'$ be its neighbour that is a leaf. Observe that for any minimum dominating set $D$ of $T$ we have $D\cap\{u,v,w,x,y,x'\}=\{v,y\}$. In particular, every minimum dominating set of $T$ is of the form $D' \cup \{v\}$ where $D'$ is a minimum dominating set of $T-\{u,v,w\}$ of size $\gamma-1$. We thus have $\Gamma(T) \leq f(\gamma-1,s) \leq \alpha^s \beta^{\gamma-s-1} \leq \alpha^s \beta^{\gamma-s}$, which gives a contradiction.
\end{proof}

\begin{Claim}\label{y_has_degree_2}
Every $4$-terminal vertex $y$ has degree $2$.
\end{Claim}

\begin{proof}
Suppose that $y$ has a neighbour $x'$ that is $k$-terminal for $k \leq 3$, different from $x$. By \Cref{y_not_support_vertex}, we have $k \geq 1$. If $k=1$, let $w'$ be a $0$-terminal neighbour of $x'$ (see \Cref{fig:path_y_one_child}). Then a set $D$ is a minimum dominating set of $T$ if and only if one of the following holds:

\begin{figure}[htbp]\centering
    			\includegraphics[height=4cm]{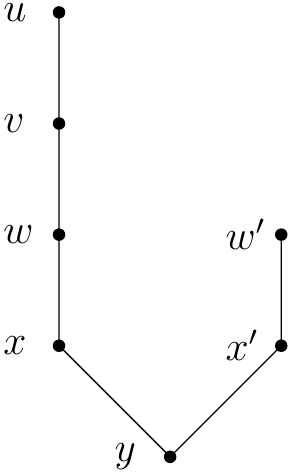}
    			\caption{Proof of \Cref{y_has_degree_2}. Case $k=1$.}
\label{fig:path_y_one_child}
	\end{figure}

\begin{itemize}
    \item either $D$ is of the form $ D_1 \cup \{ u,w,w' \}$, $D_1 \cup \{ u,w,x' \} $, $D_1 \cup \{ v,w,w' \} $ or $D_1 \cup \{ v,w,x' \}$ where $D_1$ is a minimum dominating set of $T - \{ u,v,w,x,y,w',x' \}$ of size $\gamma - 3$,
    
    \item or $D$ is of the form $D_1 \cup \{ u,x, x' \}$ or $D_2 \cup \{ u,x, w' \}$ where $D_1$ is a minimum dominating set of $T - \{ u,v,w,x,y,w',x' \}$ of size $\gamma - 3$,
    
    \item or $D$ is of the form $D_2 \cup \{ v \}$ where $D_2$ is a minimum dominating set of $T - \{ u,v,w \}$ of size $\gamma - 1$.
\end{itemize}

This gives:

\begin{align*}
    \Gamma(T)  &\leq 4 f(\gamma-3,s) + 2 f(\gamma-3,s) + f(\gamma-1,s)  \\
    &\leq \alpha^{s}\beta^{\gamma -s}(6 \beta^{-3} + \beta^{-1} )  \\
    &\leq \alpha^{s}\beta^{\gamma -s}, \\
\end{align*}
which gives a contradiction. \\

If $k=2$, let $w'$ be a $1$-terminal neighbour of $x'$, and $v'$ be a $0$-terminal neighbour of $w'$. By \Cref{Only_forks_or_3paths}, $x'$ has at most $3$ neighbours and all but one branches from $x'$ must be paths of length $2$. Let $m$ be the number of $1$-terminal neighbours of $x'$. We claim that $x'$ cannot be in a minimum dominating set. Indeed, let $D$ be a minimum dominating set of $T$, suppose $x' \in D$ and let $B$ be the union of $\{y\}$ with all the vertices of the branches from $y$ except, it it exists, the branch that contains a vertex that is $l$-terminal for some $l \geq 4$. Note that every path of length $2$ appended to $x'$ must intersect $D$ and that $|D \cap \{ y,x,w,v,u \}| \geq 2$. Hence, $|D \cap B| \geq m+3$. But then $D'=D \setminus \{u,w,x,x'\}  \cup \{v,y\}$ dominates $T$ and satisfies $|D'| \leq |D| -1$, which contradicts the definition of $D$. Hence, $x'$ cannot be in $D$, and distinguishing whether $v'$ or $w'$ is in $D$, we get that $ \Gamma(T) \leq 2 f(\gamma-1,s) \leq \alpha^{s}\beta^{\gamma -s}$, which yields a contradiction. \\

Finally, if $k=3$, then let $w'$ be a $2$-terminal neighbour of $x'$, $v'$ be a $1$-terminal neighbour of $w'$ and $u'$ be a $0$-terminal neighbour of $v'$. By \Cref{claim:strong-support}, \Cref{2_term_has_deg_2} and \Cref{path_x_deg_2}, $v'$, $w'$ and $x'$ have degree $2$. Now, any minimum dominating set must contain $v$, $v'$ and $y$, which gives $\Gamma(T) \leq f(\gamma -2, s+1) \leq \alpha^{s+1}\beta^{\gamma-3 -s} \leq \alpha^{s}\beta^{\gamma -s}$, once again which gives a contradiction.
\end{proof}

\begin{Claim}\label{Exist_4_terminal}
There exists a $4$-terminal vertex.
\end{Claim}

\begin{proof}
Suppose for a contradiction that there does not exist a $4$-terminal vertex. As $\gamma>1$, there exists a $2$-terminal vertex. If there exists a $3$-terminal vertex, then by \Cref{Only_forks_or_3paths} and \Cref{path_x_deg_2}, $T$ must be one of the graphs in \Cref{fig:3term}.

\begin{figure}[htbp]\centering
    			\includegraphics[height=2cm]{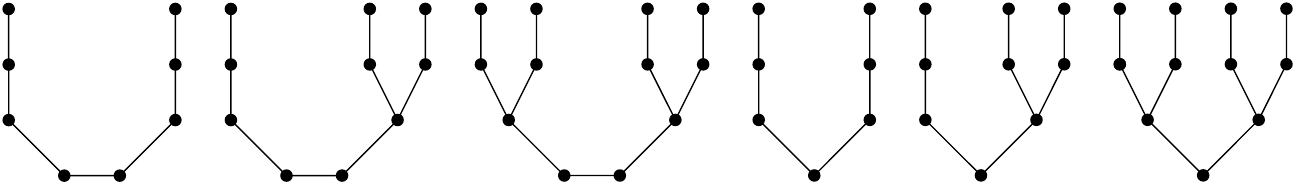}
    			\caption{All possible trees with no $4$-terminal vertex that have a $3$-terminal vertex and satisfy \Cref{Only_forks_or_3paths} and \Cref{path_x_deg_2}.}
\label{fig:3term}
	\end{figure}

One can easily check that each of them has at most $\sqrt{5}^{\gamma}$ minimum dominating sets, which gives a contradiction. \\

If there is no $3$-terminal vertex, then by \Cref{Only_forks_or_3paths}, $G$ must be one of the graphs in \Cref{fig:2term}.

\begin{figure}[htbp]\centering
    			\includegraphics[height=1.5cm]{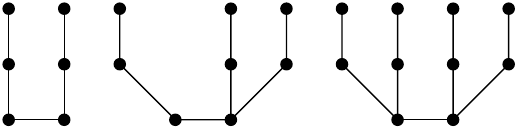}
    		    \caption{All possible trees the have no $3$-terminal vertex that have a $2$-terminal vertex and satisfy \Cref{Only_forks_or_3paths}.}
\label{fig:2term}
	\end{figure}

Once again, one can easily check that each of them has at most $\sqrt{5}^{\gamma}$ minimum dominating sets, which gives a contradiction.
\end{proof}

With the above claims, we are now able to prove \Cref{th:UB}.

\begin{proof}[Proof of \Cref{th:UB}.]

By \Cref{Exist_4_terminal}, there exists a $4$-terminal vertex $y$. Let $x$ be a $3$-terminal neighbour of $y$, $z$ be another neighbour of $y$, $w$ be a $2$-terminal neighbour of $x$, $v$ be a $1$-terminal neighbour of $w$ and $u$ be a leaf-neighbour of $v$. Note that by Claims \ref{y_has_degree_2}, \ref{path_x_deg_2}, \ref{2_term_has_deg_2} and \ref{claim:strong-support} the vertices $y$, $x$, $w$ and $v$ have no further neighbours.

First, suppose that $z$ is not a strong support vertex. Then a set $D$ is a minimum dominating set of $T$ if and only if one of the following holds:

\begin{itemize}
    \item either $D$ is of the form $D_1 \cup \{u,w\}$ or $D_1 \cup \{v,w\}$ where $D_1$ is a minimum dominating set of $T'$ of size $\gamma - 2$ where $T'$ is the tree obtained from $T$ by deleting the set of vertices $ \{ u,v,w,x,y \}$ and adding two leaves to $z$,
    \item or $D$ is of the form $D_2 \cup \{u,x\}$ where $D_2$ is a minimum dominating set of $T - \{ u,v,w,x,y \}$ of size $\gamma - 2$,
    \item or $D$ is of the form $D_3 \cup \{v\}$ where $D_3$ is a minimum dominating set of $T - \{ u,v,w \}$ of size $\gamma - 1$.
\end{itemize}

This gives:

\begin{align*}
    \Gamma(T)  &\leq 2 f(\gamma-2,s+1) + f(\gamma-2,s)+ f(\gamma-1,s)  \\
    &\leq \alpha^{s}\beta^{\gamma -s}(2 \beta^{-2}(\beta -1)^{-1}  + \beta^{-2} +  \beta^{-1} )  \\
    &\leq \alpha^{s}\beta^{\gamma -s}, \\
\end{align*}
which yields a contradiction. \\

Finally, suppose $z$ is a strong support vertex. Then a minimum dominating set $D$ of $T$ must be of the form $D_1 \cup \{ y,v\}$, or $D_1 \cup \{ x,v\}$, or $D_1 \cup \{ x,u\}$, or $D_1 \cup \{ w,v\}$, or $D_1 \cup \{ w,u\}$, where $D_1$ is a minimum dominating set of $T - \{ u,v,w,x,y \}$ of size $\gamma - 2$. This gives $f(\gamma,s) \leq 5f(\gamma-2,s) \leq \alpha^{s}\beta^{\gamma -s}$, giving again a contradiction.
\end{proof}

\subsection{Proof of the lower bound}

\Cref{LowerBound} is an immediate corollary of the following theorem.

\begin{theorem}
For each non-negative integer $k$ one can find a tree with domination number $\gamma=2k+1$ that has $5^k+3^k$ minimum dominating sets and a tree with domination number $\gamma=2k+2$ with $2\cdot 5^k + 3^k$ minimum dominating sets.
\end{theorem}

\begin{figure}[htbp]\centering
    			\includegraphics[height=6cm]{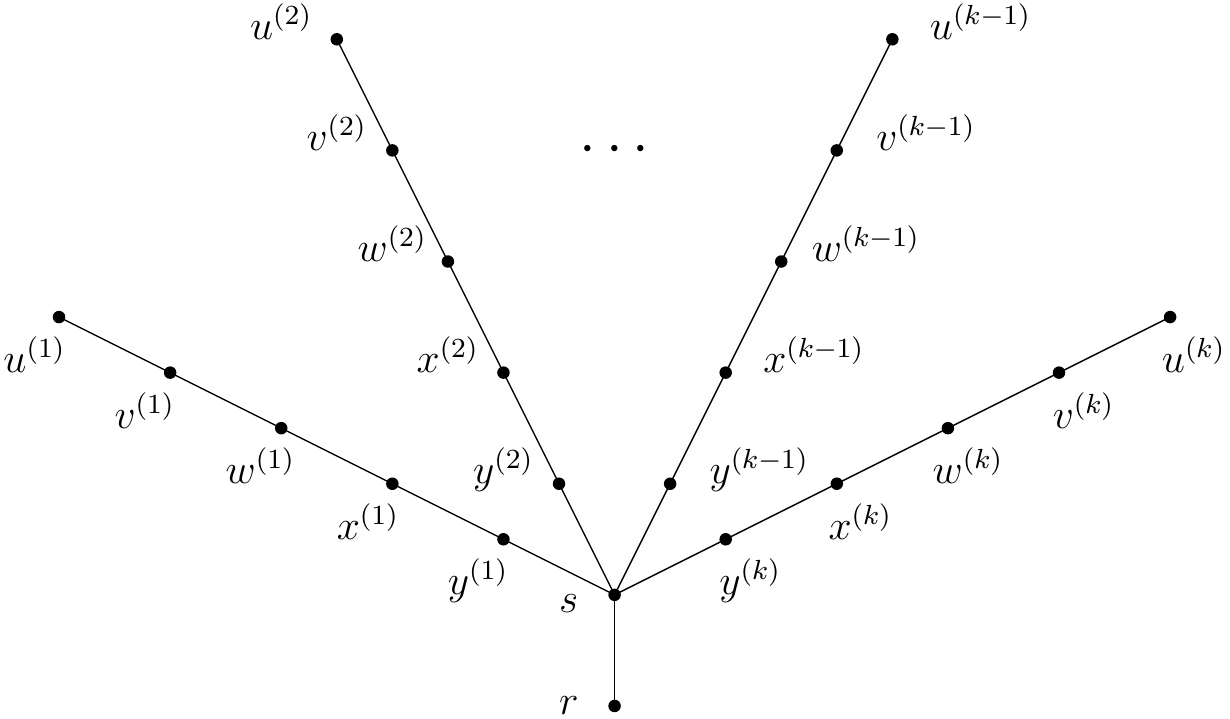}
    			\caption{The tree $G_k$ with domination number $\gamma=2k+1$ and $5^k+3^k$ minimum dominating sets. The tree $H_k$ with domination number $\gamma=2k+2$ and $2 \cdot 5^k+3^k$ minimum dominating set is formed from $G_k$ by appending a path on three vertices to $s$.}
\label{fig:LB}
	\end{figure}

\begin{proof}
Let $G_k$ be a graph on vertices $V=\{r,s\}\cup\bigcup_{i=1}^kV^{(i)}$, where $V^{(i)}=\{y^{(i)},x^{(i)},w^{(i)},v^{(i)},u^{(i)}\}$, and with edges $E=\{rs\}\cup\bigcup_{i=1}^k\{sy^{(i)},y^{(i)}x^{(i)},x^{(i)}w^{(i)},w^{(i)}v^{(i)},v^{(i)}u^{(i)}\}$. See \Cref{fig:LB} for an illustration.

Let $D$ be a minimum dominating set of $G_k$. It must satisfy $|D \cap \{r,s\}|\geq 1$ and also $|D \cap V^{(i)}|\geq 2$ for all $i$, hence $|D|\geq 2k+1$. On the other hand, there clearly exist dominating sets of $G_k$ of size $2k+1$, for example $\{s\}\cup\bigcup_{i=1}^k\{x^{(i)}, v^{(i)}\}$. Therefore, the domination number of $G_k$ is $2k+1$ and for each minimum dominating set $D$ of $G_k$ we have $|D \cap \{r,s\}|= 1$ and $|D \cap V^{(i)}|= 2$ for all $i$.

If $D \cap \{r,s\}=s$, then for each $i$, $D \cap V^{(i)}$ is one of $\{y^{(i)},v^{(i)}\}$, $\{x^{(i)},v^{(i)}\}$, $\{x^{(i)},u^{(i)}\}$, $\{w^{(i)},v^{(i)}\}$, $\{w^{(i)},u^{(i)}\}$.
If, on the other hand, $D \cap \{r,s\}=r$, then for each $i$, $D \cap V^{(i)}$ is one of $\{y^{(i)},v^{(i)}\}, \{x^{(i)},v^{(i)}\}, \{x^{(i)},u^{(i)}\}$. Hence $G_k$ has $5^k+3^k$ minimum dominating sets.

For the second part, let $H_k$ be the graph on vertices $V(G_k)\cup\{u,v,w\}$ and with edges $E(G_k)\cup\{uv,vw,ws\}$. Because of $u$, for each dominating set of $H_k$ we have $|D\cap\{u,v,w\}|\geq1$. Also, $\{s,v\}\cup\bigcup_{i=1}^k\{x^{(i)},v^{(i)}\}$ is a dominating set of $H_k$ of size $2k+2$. By similar arguments to the case of $G_k$ we obtain that the domination number of $H_k$ is $2k+2$ and that the number of minimum dominating sets is $2\cdot5^k+3^k$.
\end{proof}

\section{Minimum total dominating sets in forests}

Let $G(k,l,m)$ be the following graph: we start with a vertex $y$ and we append to it vertices $x_1, \dots ,x_k$. Then, to each $x_a$ we append vertices $w_{a,1}, \dots, w_{a,k}$ and each $w_{a,b}$ gets itself a $1$-terminal neighbour $v_{a,b}$ with $m$ leaves $u_{a,b,1}, \dots ,u_{a,b,m}$ (see \Cref{fig:G_{klm}}). We will prove that for a suitable choice of $k$, $l$ and $m$, this graph is a counterexample to \Cref{Conj_Total_UB}.

\begin{figure}[htbp]\centering
    			\includegraphics[height=4.5cm]{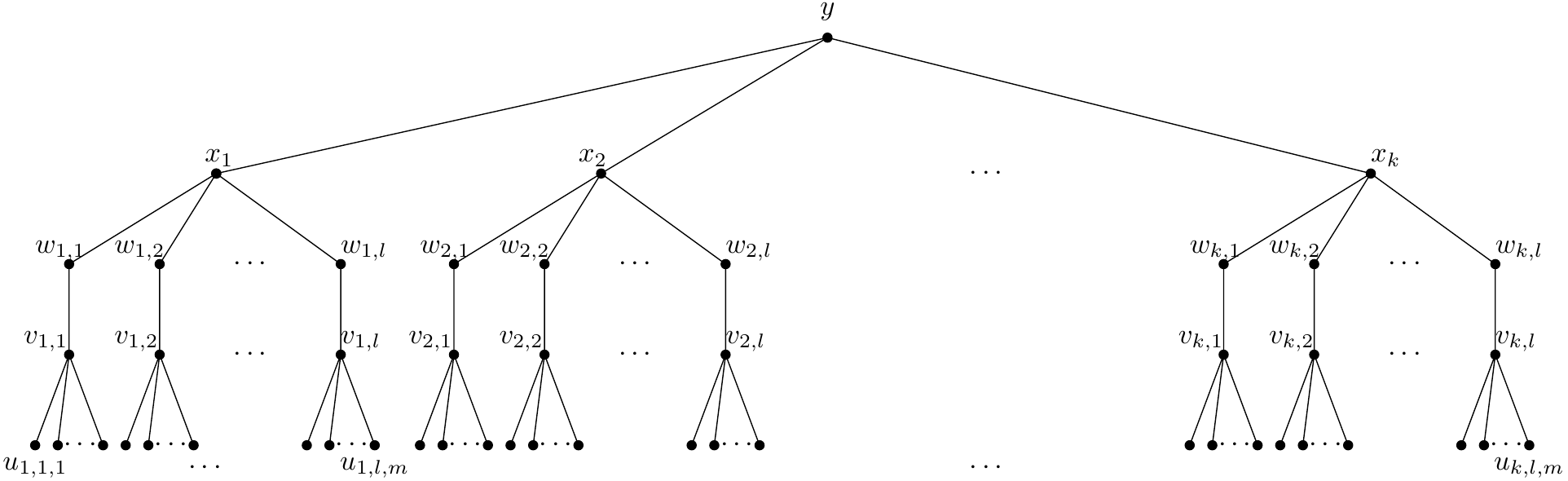}
    			\caption{The tree $G(k,l,m)$.}
\label{fig:G_{klm}}
	\end{figure}

\begin{lemma}\label{CounterExampleHenning}
The tree $G=G(k,l,m)$ satisfies $n(G)=klm+2kl+k+1$, $\gamma_t(G) =2kl+1 $, $\Gamma_t(G) = k( (m+1)^l - m^l)^k$.
\end{lemma}

\begin{proof}
The equality $n(G)=klm+2kl+k+1$ is trivial. For the second equality, let $D$ be a minimum total dominating set of $G(k,l,m)$. Each of the $v_{a,b}$, $1 \leq a \leq k$ and $1 \leq b \leq l$, must be in $D$. Moreover, $y$ and each of the $v_{a,b}$ must be dominated, hence $D$ must contain an element of $N(y)$ and of each of the $N(v_{a,b})$. Those sets are pairwise disjoint and do not contain any of the $v_{a,b}$. Hence, $\gamma_t = |D| \geq 2kl+1$. Moreover, $D= \{ w_{a,b} :  1 \leq a \leq k, 1 \leq b \leq l\} \cup \{ v_{a,b} :  1 \leq a \leq k, 1 \leq b \leq l\} \cup \{ x_1 \}$ is clearly a total dominating set, which proves $\gamma_t(G) =2kl+1$.

For the third equality, let $D$ be a minimum total dominating set of $G(k,l,m)$. As stated in the paragraph above, $D$ must contain the $kl$ vertices $v_{a,b}$, an element of $N(y)$ and an element of each of the $N(v_{a,b})$. Each of the $x_i$ must also be dominated, so for each $1 \leq i \leq k$, one of the $w_{a,i}$ for $1 \leq a \leq l$ must be in $D$. Conversely, it is easy to see that a set $D$ satisfying the above conditions is a minimum total dominating set of $G$. Consequently, a counting argument shows that $\Gamma_t(G) = k ( (m+1)^l-m^l)^k$.
\end{proof}

\begin{lemma}\label{Disprove_Conj_Henning}
The tree $T=G(3,16,3)$ satisfies $\Gamma_t(T) \geq c^{\gamma_t} \left(\frac{n-\frac{\gamma_t}{2}}{\frac{\gamma_t}{2}} \right)^\frac{\gamma_t}{2}$ for some $c>1$.
\end{lemma}

\begin{proof}
By \Cref{CounterExampleHenning}, for $k=m=3$ and $l=16$, we get $\Gamma_t(T) \approx 2.306 \times 10^{29}$, whereas $\left(\frac{n-\frac{\gamma_t}{2}}{\frac{\gamma_t}{2}} \right)^\frac{\gamma_t}{2} \approx 2.302 \times 10^{29}$, whence $c={\left(\Gamma_t(T)/\left(\frac{n-\frac{\gamma_t}{2}}{\frac{\gamma_t}{2}} \right)^\frac{\gamma_t}{2}\right)}^{1/97}$ is greater than $1$. 
\end{proof}

We note that $c > 1.000017$ and that \Cref{Disprove_Conj_Henning} already disproves \Cref{Conj_Total_UB}. \Cref{th:total_LB} follows:

\begin{proof}[Proof of \Cref{th:total_LB}.]
Consider the tree $G_N$ made of $N$ copies $G^{(1)}, \dots G^{(N)}$ of the tree $T=G(3,16,3)$ where we put an edge between an arbitrary leaf $G^{(i)}$ (say $u_{3,16,3}^{(i)}$) and an arbitrary leaf of $G^{({i+1)}}$ (say $u_{1,1,1}^{(i+1)}$) for all $1 \leq i \leq N-1$. Each total dominating set of $G$ contains all $v_{a,b}^{(i)}$ as they are support vertices. Therefore, in each total dominating set of $G$, $u_{1,1,1}^{(i)}$ and $u_{3,16,3}^{(i)}$ are totally dominated by $v_{1,1}^{(i)}$ and $v_{3,16}^{(i)}$, respectively. Thus, a set $D \subset V(G)$ is a total dominating set of $G$ if and only if for all $i=1,...,N$ the set $D \cap V(G^{(i)})$ is a total dominating set of $G^{(i)}$. From this we have $\gamma_t=\gamma_t(G_N)=N\gamma_t(T)$ and $\Gamma_t(G_N)=\Gamma_t(T)^N$. In view of $n=n(G_N)=Nn(T)$ and \Cref{Disprove_Conj_Henning}, we see that $G_N$ is a tree satisfying $\Gamma_t(G_N) \geq c^{\gamma_t} \left(\frac{n-\frac{\gamma_t}{2}}{\frac{\gamma_t}{2}} \right)^\frac{\gamma_t}{2}$, where $c>1$ is the same constant as in \Cref{Disprove_Conj_Henning}.
\end{proof}

\section{Concluding remarks}

In this section, we lay out two questions we deem interesting for further research. Consider for a forest $F$ with total domination number $\gamma_t$ the quantity
\begin{align*}
    \lambda(F) = \Gamma_t(F)^{\frac{1}{\gamma_t}}
    \left(\frac{n-\frac{\gamma_t}{2}}{\frac{\gamma_t}{2}}\right)^{-\frac{1}{2}}.
\end{align*}

Whenever $\lambda(F)>1$, $F$ fails \Cref{Conj_Total_UB}, as is the case for $F=G(3,16,3)$ where $\lambda(F) > 1.000017$. On the other hand, \Cref{th:total_ub} by Henning et al. asserts that $\lambda(F)\leq 8\sqrt{e}$ for any forest, so it is natural to ask what the supremum of $\lambda(F)$ over all forests is. 

Considering that the $G(k,l,m)$ trees have a more complex structure than the trees $G_k$ that show the lower bound on the number of minimum dominating sets, it seems that determining the supremum of $\lambda$ could be hard. However, there are some natural first steps to take. First of all, there is good reason to believe that the supremum should not change if we take it over trees only, since we may connect leaves of strong support vertices of different components without changing either the number of vertices or the total domination number. Furthermore, we believe the $G(k,l,m)$ trees may in some sense be the worst that can happen to $\lambda$. Indeed, when attempting to bound $\Gamma_t(T)$ by recursively pruning parts of the tree as for the dominating sets, one naturally arrives at such structures as an obstacle. However, there is no reason to assume that there should be a single $G(k,l,m)$-tree that achieves the supremum. This leads us to the following question.

\begin{Question}
Let $\mathcal{F}$ be the set of all forests and $\mathcal{F'}$ be the set of forests in which each tree is a $G(k,l,m)$-tree. Is it true that
\begin{align*}
    \sup_{F\in \mathcal{F}} \lambda(F) = \sup_{F\in \mathcal{F'}} \lambda(F)?
\end{align*}
\end{Question}

Since the above question is borne out of the suspicion that, in a loose sense, the $G(k,l,m)$ trees maximize $\lambda$, it is also natural to ask for a stability result manifesting this suspicion. We cannot hope that every $F$ with $\lambda(F)$ above a certain threshold belongs to $\mathcal{F'}$, since it is always possible to add trees with a small number of total dominating sets to a large forest $F\in \mathcal{F'}$ without decreasing $\lambda(F)$ too much. The following question is formulated with this degeneracy in mind.

\begin{Question}
Does there exist for every $\varepsilon>0$ a constant $c> 1$ such that every forest $F$ with $\lambda(F)> c$ contains a subforest $F'\in \mathcal{F'}$ with $\lambda(F')\geq \lambda(F)$ and $\vert V(F') \vert \geq (1-\varepsilon)\vert V(F) \vert$?
\end{Question}

\section*{Acknowledgement}

The authors would like to thank Professor Béla Bollobás for his valuable comments.

\bibliographystyle{abbrvnat}  
\bibliography{bib}

\Addresses

\end{document}